\documentclass[12pt]{amsart}
\usepackage{amsfonts}
\usepackage{amssymb}
\newtheorem{thm}{Theorem}[section]
\newtheorem{cor}[thm]{Corollary}
\newtheorem{lem}[thm]{Lemma}
\newtheorem{prop}[thm]{Proposition}
\theoremstyle{definition}
\newtheorem{defn}[thm]{Definition}
\theoremstyle{remark}
\newtheorem{rem}[thm]{Remark}

\begin{document}

\title[Directed Riemannian manifolds ]
{Directed Riemannian manifolds of pointwise constant relative
sectional curvature}%

\author{ Georgi Ganchev and Vesselka Mihova}%

\address{Bulgarian Academy of Sciences, Institute of Mathematics and Informatics,
Acad. G. Bonchev Str. bl. 8, 1113 Sofia, Bulgaria}%
\email{ganchev@math.bas.bg}
\address{Faculty of Mathematics and Informatics, University of Sofia,
J. Bouchier Str. 5, 1164 Sofia, Bulgaria}
\email{mihova@fmi.uni-sofia.bg}%




\begin{abstract}
We study a class of Riemannian manifolds with respect to the
covariant derivative of their curvature tensors. We introduce
geometrically the class of directed Riemannian manifolds of
pointwise constant relative sectional curvature and give a tensor
characterization for such manifolds. We prove that all rotational
hypersurfaces are directed and find the rotational hypersurfaces of
pointwise constant relative sectional curvature. For the class of
directed Riemannian manifolds of pointwise constant relative
sectional curvature having a totally umbilical scalar distribution we
prove a structural theorem and a theorem of Schur's type.
\end{abstract}

\maketitle
\thispagestyle{empty}

\section{Introduction}

Let $(M, g)$ be a Riemannian manifold with Levi-Civita connection $\nabla$
and curvature tensor $R$. If $V$ is an $n$-dimensional vector space identified
with the tangent space at an arbitrary point in $M$, denote by ${\mathcal R}(V)$
the linear space of all tensors of type (0,4) over $V$ having the symmetries
of $R$.

According to the general theory of group representations \cite{W} there
exists a splitting of ${\mathcal R}(V)$ into irreducible components
under the action of $O(n)$. Singer and Thorpe  \cite{ST} and Nomizu
\cite{N} give explicitly a decomposition of ${\mathcal R}(V)$ and
describe it geometrically in terms of the well-known classes of
Riemannian manifolds of constant sectional curvature, Einstein
manifolds and conformally flat Riemannian manifolds.

In the case when $(M,g,J)$ is an almost Hermitian manifold with almost complex structure
$J$ and $V$ is a $2n$-dimensional Hermitian  vector space, Tricerri and
Vanhecke   give in \cite{TV} a complete explicit decomposition of ${\mathcal R}(V)$
under the action of $U(n)$. In this case the splitting of ${\mathcal R}(V)$
gives many new classes of almost Hermitian manifolds with respect to $R$ and
leads to the problem of their geometrical description.

Following this scheme of studying Riemannian manifolds it seems natural to
investigate the linear space $\nabla {\mathcal R}(V)$ of all tensors of type (0,5)
over $V$ having the symmetries of the covariant derivative $\nabla R$ of the
curvature tensor $R$ of a Riemannian manifold $(M,g)$. A complete explicit
decomposition of $\nabla {\mathcal R}(V)$ under the action of $O(n)$ has been
given by Gray and Vanhecke in \cite{GV}. The zero space of this splitting
$(\nabla R = 0)$ leads to the class of locally symmetric Riemannian manifolds
and this class corresponds to the class of locally flat Riemannian manifolds
which is the zero class $(R =0)$ in the splitting of ${\mathcal R}(V)$. However
we have to mention that for the classes of Riemannian manifolds with
respect to $\nabla R$ $(\nabla R \not = 0)$ it is not known very much.

Conformally flat Riemannian manifolds ($M, g, d\tau$) with metric $g$ and
scalar 1-form $d\tau$ ($\tau$ being the scalar curvature
of $(M,g)$) have been studied in \cite{GM}.

In this paper we consider the class of Riemannian manifolds whose covariant
derivative $\nabla R$ of the curvature tensor is  constructed only by the
metric $g$ and the scalar 1-form $d\tau$ . This class corresponds
to the class of Riemannian manifolds of
constant sectional curvature (in the splitting of ${\mathcal R}(V)$).
We introduce geometrically the class of directed Riemannian manifolds of
pointwise constant relative sectional curvature and prove that these manifolds
form the class of Riemannian manifolds with special covariant derivative
$\nabla R$ of the curvature tensor mentioned above. We prove that any
rotational hypersurface is a directed Riemannian manifold and find all
rotational hypersurfaces of pointwise constant relative sectional curvature.
For the special subclass of the directed Riemannian manifolds of pointwise
constant relative sectional curvature whose distribution is totally umbilical
we prove a structural theorem and a theorem of Schur's type.

\section{Directed Riemannian manifolds of pointwise constant relative
sectional curvature}

Let $(M,g)$ be a Riemannian manifold with Levi-Civita connection
$\nabla$. The Riemannian curvature operator $R$ is given by $R(X,Y)
= [\nabla_X , \nabla_Y ] - \nabla _{[X,Y]}$ and the corresponding
curvature tensor of type (0,4) is defined by $R(X,Y,Z,U) =
g(R(X,Y)Z,U)$ for arbitrary differentiable vector fields X,Y,Z,U.
Further the algebra of all differentiable vector fields on $M$ will
be denoted by ${\mathcal X}M$.

The covariant derivative $\nabla R$ of the curvature tensor $R$ has the
following symmetries:
$$\begin{array}{l}
(\nabla _{W}R)(X,Y,Z,U) = - (\nabla _{W}R)(Y,X,Z,U) = - (\nabla
_{W}R)(X,Y,U,Z);\\
[2mm]
\sigma _{XYZ} (\nabla _{W}R)(X,Y,Z,U) = 0;\\
[2mm]
\sigma _{WXY} (\nabla _{W}R)(X,Y,Z,U) = 0,
\end{array}\leqno(1)$$
where $W,X,Y,Z,U \in {\mathcal X}M$ and $\sigma$ denotes the corresponding cyclic
summation.

We denote by $\tau$ the scalar curvature of the manifold $(M,g)$ and by $\pi$ the
tensor
$$\pi (X,Y,Z,U) = g(Y,Z)g(X,U) - g(X,Z)g(Y,U); \hspace{0,5cm} X,Y,Z,U \in
{\mathcal X}M.$$

We recall that a Riemannian manifold of constant
sectional curvatures is characterized by the equality
$$ R = {\frac{\tau }{n(n-1)}}\pi, \leqno{(2)}$$
i.e. the curvature tensor of a Riemannian manifold of constant sectional
curvatures is constructed only by the metric $g$.

Let $\omega$ be a 1-form on the Riemannian manifold $(M, g)$. We consider the tensor
$$\begin{array}{l}
\Pi (\omega )(W,X,Y,Z,U) = 2\omega (W)\pi (X,Y,Z,U) + \omega (X)\pi (W,Y,Z,U)\\
[2mm]
+ \omega (Y)\pi (X,W,Z,U) + \omega (Z)\pi (X,Y,W,U) + \omega (U)\pi (X,Y,Z,W).
\end{array}$$
It is easy to check that the tensor $\Pi (\omega )$ has the symmetries (1) of
the tensor $\nabla R$.

Our aim in this paper is to study the class of Riemannian manifolds characterized by
the condition
$$ \nabla R = {\frac{1}{2(n-1)(n+2)}}\Pi (\omega ). \leqno{(3)}$$

With respect to $\nabla R$ this class formally corresponds to the class of
Riemannian manifolds of constant sectional curvatures.
In terms of the decomposition of $\nabla R$ \cite{GV} the condition (3) means that
$\nabla R$ coincides with its component in the space $\nabla {\mathcal R}_{I}$.

In this section we characterize  the equality (3) geometrically.

Let $E = span\{X,Y\}$ be a 2-plane in the tangent space $T_{p}M$ at a point $p$
in $M$ and $\{X,Y\}$ be an orthonormal basis of $E$. The tensor $\nabla R$
generates the  1-form $\varphi _{E}$ defined on $E$ as follows:
$$ \varphi _{E}(Z) = (\nabla _{Z}R)(X,Y,Y,X), \hspace{0,5cm} Z \in E. \leqno{(4)}$$

Because of the properties (1) of $\nabla R$ the 1-form $\varphi _{E}$ does not
depend on the orthonormal basis of $E$.

The 1-form $\varphi _{E}$ defined on $E$ by (4) is said to be {\it a sectional 1-form}.

Let now $\eta $ be a unit 1-form on $(M,g)$ and $\Delta $ be the
distribution of $\eta $, i.e.
$$ \Delta (p) = \{ X \in T_{p}M : \eta (X) = 0 \},\; p \in M.$$

\begin{defn} A Riemannian manifold $(M,g)$ is said to be {\it directed}
if there exists a unit 1-form $\eta $ on $M$ such that

i) $\varphi _{E} = k(E,p)\,\eta \vert _{E}\;$ for all $\,E \not \subset \Delta;$

ii) $\varphi _{E} = 0\;$ for all $\,E \subset \Delta$.

 For any 2-plane $E \not \subset \Delta$ the function
 $k(E,p)$ is said to be {\it a relative sectional curvature}.
\end{defn}

The condition i) means that all sectional 1-forms
$\varphi _{E},\; E \not \subset \Delta\,$ are collinear with the restriction
of the 1-form $\eta $ to the 2-plane $E$. We say that $(M,g)$ {\it is directed
by the 1-form $\eta $}.

\begin{defn} A directed Riemannian manifold $(M,g)$ is said to be
{\it of pointwise constant relative sectional curvatures} if the relative
sectional curvature $k(E,p)$ of any 2-plane $E \not \subset \Delta$ does not
depend on $E$.
\end{defn}

In order to find a tensor characterization for the Riemannian manifolds
described in Definition 2.2 we need the following
\begin{lem}\label{lem 1}
Let $L$ be a tensor of type $(0,5)$ satisfying the following equalities
\begin{itemize}
\item[i)] $L(W,X,Y,Z,U) = -L(W,Y,X,Z,U) = -L(W,X,Y,U,Z);$
\vskip 2mm

\item[ii)] $\sigma _{XYZ}L(W,X,Y,Z,U) = 0;$
\vskip 2mm

\item[iii)] $\sigma _{WXY}L(W,X,Y,Z,U) = 0$
\end{itemize}
for all $W,X,Y,Z,U \in {\mathcal X}M$.

If $L(X,X,Z,Z,X) = 0$ for arbitrary $X,Z \in {\mathcal X}M$, then $L \equiv 0.$
\end{lem}

{\it Proof}. Substituting successively $X$ by $X + Y$ and by $X - Y$ into the equality
$$L(X,X,Z,Z,X) = 0$$
and taking into account the properties of $L$, we obtain
$$ L(X,Y,Z,Z,Y) + 2L(Y,X,Z,Z,Y) = 0. \leqno{(5)}$$
This implies that
$$ L(Y,X,Z,Z,Y) = L(Z,X,Y,Y,Z). \leqno{(6)}$$

Applying the condition iii) to $L(X,Y,Z,Z,Y)$ and taking into account
(6) we find
$$ L(X,Y,Z,Z,Y) - 2L(Y,X,Z,Z,Y) = 0.$$

The last equality combined with (5) implies $L(X,Y,Z,Z,Y) = 0$ for all
$X,Y,Z \in {\mathcal X}M$. Now it follows in a standard way that $L \equiv 0$.
\hfill {$\square$}

We give a tensor characterization for directed manifolds of pointwise constant
relative sectional curvatures.
\begin{thm}
Let $(M,g)$ be a Riemannian manifold directed by the unit 1-form $\eta $. Then
$(M,g)$ is of pointwise constant relative sectional curvatures $k(p)$ if and only if
$$ \nabla R = {\frac{1}{4}}k(p)\Pi (\eta ).\leqno{(7)}$$
The function $k(p)$ satisfies the equality
$$ d\tau = {\frac{(n-1)(n+2)}{2}}\, k\,\eta, \leqno{(8)}$$
where $\tau $ is the scalar curvature of the manifold.
\end{thm}

{\it Proof.} To prove the first implication we put
$$ L = \nabla R - {\frac{1}{4}}\,k\,\Pi (\eta ).$$

Under the conditions of the theorem it is easy to check that
$L (X,X,Y,Y,X) = 0$ for all $X,Y \in {\mathcal X}M$. Applying Lemma 2.3 we obtain (7).

The inverse is an easy verification.

The equality (8) follows from (7) by two contractions.

\hfill {$\square$}

Theorem 2.4 implies immediately
\begin{cor}\label{cor 1}
Let $(M,g)$ be a directed Riemannian manifold of pointwise constant relative
sectional curvatures. Then, $(M,g)$ is locally symmetric if and only if $d\tau = 0$.
\end{cor}

Considering directed Riemannian manifolds of pointwise constant relative
sectional curvatures $k$ and of nonconstant scalar curvature $\tau$, i. e.
$d\tau \not = 0$ on $M$, we compute (up to an orientation of $\eta $) from  (8)
$$ k = {\frac{2\,\Vert d\tau \Vert}{(n-1)(n+2)}},\quad
\eta = \frac{1}{||d\tau||}\,d \tau\,. \leqno{(9)}$$

Hence, the 1-form $\eta $ is uniquely determined by the metric $g$.

\section{Examples}
In this section we give examples of the manifolds introduced in the previous
section among the rotational hypersurfaces.

First we need some formulas.

Let $(M,g)$ be a rotational hypersurface in the Euclidean space
${\bf R}^{n+1}$ with a rotational axis oriented by a unit vector $e$. We consider $M$
as a 1-parameter family of spheres $S^{n-1} (t), t \in J$, given by the
equalities
$$ (X - x_0(t))^{2} = r^{2}(t), \hspace{0,5cm} e(X - x_0(t)) = 0,\leqno{(10)}$$
where $x_0(t)$ and $r(t)$ are respectively centers and radii of the spheres.
Further we assume that the rotational hypersurface $M$ is also given by a
vector-valued function $X(u^{1},...,u^{n-1},t)$ satisfying (10), where
${u^{1},...,u^{n-1},t}$ is a local coordinate system on $M$.

Taking partial derivatives of (10) we find
$$ (X - x_0)X_{\alpha } = 0, \hspace{5mm} eX_{\alpha } = 0; \alpha = 1,...,n-1,$$
$$(X - x_0)X_{t} = rr', \hspace{0,5cm} eX_{t} = 1.$$
Then the vector $X - x_0 - rr'e \,$ is normal to $M$ at the point $X$ and we can
choose the unit normal to $M$ by the equality
$$ N = - {\frac{X - x_0 - rr'.e}{r\sqrt {1+r'^{2}}}}.$$

Denote by $\xi $ the unit vector field tangent to $M$ and perpendicular to
the parallels $S^{n-1} (t)$. Up to a sign we have
$$ \xi = \sqrt {1+r'^{2}} e - r' N.$$

If $\nabla '$ is the standard flat connection in ${\bf R}^{n+1}$ we find
the Weingarten formulas on $M$:
$$ \nabla '_{x}N = {\frac{1}{r \sqrt {1+r'^2}}} x, \hspace{0,5cm} x \perp \xi;$$
$$ \nabla '_{\xi }N = {\frac{r''}{(\sqrt {1+r'^2})^3}} \xi.$$

Hence,the second fundamental tensor $h$ of $M$ has the following structure
$$ h = {\frac{1}{r \sqrt {1+r'^2}}} g - {\frac{1+r'^2+rr''}{r(\sqrt
{1+r'^2})^3}}\eta \otimes \eta, \leqno{(11)} $$
where $\eta $ is the dual 1-form of the unit vector field $\xi $.

Substituting $h$ from (11) into the Gauss equation we find the curvature tensor
of any rotational hypersurface has the following form (see also \cite{GM}):
$$ R = a\pi + b\Phi , \leqno{(12)} $$
where $a$ and $b$ are the functions
$$ a = {\frac{1}{r^{2}(1+r'^{2})}}, \hspace{0,5cm} b = - {\frac{1+r'^{2}+rr''}
{r^{2}(1+r'^{2})^{2}}} \leqno{(13)}$$
and $\Phi $ is the tensor
$$ \Phi (X,Y,Z,U) = g(Y,Z)\eta (X)\eta (U) - g(X,Z)\eta (Y)\eta (U)$$
$$ +g(X,U)\eta (Y)\eta (Z) - g(Y,U)\eta (X)\eta (Z); \hspace{0,5cm}
X,Y,Z,U \in {\mathcal X}M.$$

Let $\nabla $ be the Levi-Civita connection of the rotational hypersurface
$(M,g)$. Applying the second Bianchi identity to (12) we obtain
$$ \nabla _{x}\xi = \lambda x, \hspace{0,5cm} \lambda = {\frac{\xi (a)}
{2b}}, \hspace{0,5cm} x \perp \xi; \leqno{(14)}$$
$$ (\nabla _{X} \eta )(Y) = \lambda [g(X,Y) - \eta (X)\eta (Y)], \hspace{0,5cm}
X,Y \in {\mathcal X}M; \leqno{(15)}$$
$$ da = \xi (a).\eta = 2\lambda b.\eta; \leqno{(16)}$$
$$ db = \xi (b).\eta; \leqno{(17)}$$

Taking into account (12), (15), (16) and (17) we calculate with respect to local
coordinates
$$ \nabla _{i}R_{jkpq} = \lambda b(2\eta _{i}\pi _{jkpq} + \eta _{j}\pi _{ikpq}
 + \eta _{k}\pi _{jipq} \leqno{(18)}$$
$$ + \eta _{p}\pi _{jkiq} + \eta _{q}\pi _{jkpi})
+ (\xi (b) - 2b\lambda )\eta _{i}\Phi _{jkpq}.$$

If $E = span\{X,Y\}$ is an arbitrary 2-plane in $T_{p}M, p \in M$ with an
orthonormal basis $\{X,Y\}$, we denote by $\gamma $ the angle between $\xi $ and
$E$. Then we have
$$ \cos ^2\gamma = \eta ^2(X) + \eta ^2(Y).$$
Taking into account the defining equality (4) from (18) we obtain
$$ \varphi _{E} = [4\lambda b + (\xi (b) - 2b\lambda)\cos ^2\gamma ]\eta.$$

Thus, we have
\begin{prop}
Every rotational hypersurface is a directed Riemannian manifold.
\end{prop}

Now we shall find the rotational hypersurfaces of pointwise constant relative
sectional curvature.

As a consequence of (18), (16) and Theorem 2.4 we obtain
\begin{prop}
A rotational hypersurface $(M,g)$ with curvature tensor (12) is of pointwise
constant relative sectional curvature iff
$$ a - b = B = const $$
\end{prop}

By use of the formulas (13) we find
$$ a - b = {\frac{2(1+r'^2)+rr''}{r^2(1+r'^2)^2}} = B. \leqno{(19)}$$

Solving the differential equation (19) we obtain
\begin{prop}
A rotational hypersurface $(M,g)$ with meridian $t = t(r)$ is of pointwise
constant relative sectional curvature iff
$$ t = \int {\frac{r\sqrt{Ar^2+B}}{\sqrt{1-Ar^4-Br^2}}}dr, \hspace{0,5cm}
0 < (Ar^2+B)r^2 < 1. \leqno{(20)} $$
\end{prop}

Putting $u^2 = Ar^2+B, m = {\frac{\sqrt {B^2+4A} - B}{2A}},
m' = {\frac{\sqrt {B^2+4A} + B}{2A}},$ we obtain the meridian of the
hypersurface has equations
$$ r = \sqrt {{\frac{u^2-B}{A}}}, \hspace{0,5cm} t = {\frac{1}{A}}\int
{\frac{u^2}{\sqrt {(1-mu^2)(1+m'u^2)}}}du.$$

Further we consider the cases:

I) $A > 0.$ Putting $u = \sqrt {{\frac{1-x^2}{m}}}, x \in (0,1)$ we find the
meridian of the rotational hypersurface has the following equations:
$$ r =\sqrt {{\frac{1-x^2}{m}}-{\frac{B}{A}}}, \hspace{0,5cm}
t = {\frac{-1}{Am\sqrt {m+m'}}}(J_1 - J_2), \leqno{(21)}$$
where
$$ J_1 = \int {\frac {dx}{\sqrt {(1-x^2)(1-k^2x^2)}}}, \hspace{0,5cm}
J_2 = \int {\frac {x^2dx}{\sqrt {(1-x^2)(1-k^2x^2)}}}, \hspace{0,5cm}
(k = \sqrt {{\frac{m'}{m+m'}}} < 1)$$
are the integrals of Legendre of first type and of second type, respectively.

II) $A < 0.$ Putting $u = {\frac{x}{\sqrt {-m'}}}, x \in (0,1)$ we find the
equations
$$ r = \sqrt {-{\frac{x^2}{m'A}}-{\frac{B}{A}}}, \hspace{0,5cm}
t = -{\frac{1}{Am\sqrt {-m'}}}J_2. \leqno{(22)}$$

\section{The case of a totally umbilical distribution}
Let $(M,g)$ be a Riemannian manifold with a unit vector field $\xi $. By
$\eta $ and $\Delta $ we denote respectively the dual to $\xi $ 1-form and
the distribution perpendicular to $\xi $.

The distribution $\Delta $ is said to be {\it totally umbilical} if
$$ \nabla _{x}\xi = \lambda x, \leqno{(23)}$$
where $x \in \Delta $ and $\lambda $ is a function on $M$.

From (14) it follows that every rotational hypersurface has a totally umbilical
distribution.

If we set $\theta (X) = d\eta (\xi ,X), \hspace{0,5cm} X \in {\mathcal X}M$, then
from (23) it follows that
$$ (\nabla _{X}\eta )(Y) = \lambda [g(X,Y) - \eta (X)\eta (Y)] +
\eta (X)\theta (Y). \leqno{(24)}$$
and
$$ d\eta = \eta \wedge \theta. \leqno{(25)}$$

The last equality means that the distribution $\Delta $ is involutive.

Taking into account (24) we find the Gauss formula
for the distribution $\Delta $:
$$ \nabla _{x}y = D_{x}y - \lambda g(x,y)\,\xi; \hspace{0,5cm} x,y \in \Delta,
\leqno{(26)}$$
where $D$ is the Levi-Civita connection of the distribution $\Delta $.

Next we denote by $K$ the curvature tensor of $D$ and find the Gauss
equation for the distribution $\Delta $:
$$ R(x,y,z,u) = K(x,y,z,u) - \lambda ^2\pi (x,y,z,u); \hspace{0,5cm}
x,y,z,u \in \Delta. \leqno{(27)}$$
Taking into account (26) and (27) we calculate
$$ (\nabla _{w}R)(x,y,z,u) = (D_{w}K)(x,y,z,u) + d\lambda ^2(w)\pi (x,y,z,u);
w,x,y,z,u \in \Delta. \leqno{(28)}$$

When $d\tau \not = 0$ on the Riemannian manifold $(M,g)$ the distribution
of the 1-form $d\tau $ is said to be {\it the scalar distribution}.

Now we can prove
\begin{thm}\label{th 2}
Let $(M,g)$ be a connected directed Riemannian manifold of pointwise constant
relative sectional curvature. If the scalar distribution of the manifold is
totally umbilical, then $M$ is a one-parameter family of locally symmetric
submanifolds.
\end{thm}

{\it Proof.} Because of  (25) the distriburion $\Delta$ is involutive.

Let $p \in M$ and $S_p $ be the maximal integral submanifold of
the distribution $\Delta $ through the point $p$. Since $R$ and $K$ satisfy
the second Bianchi identity, then the equality (28) implies $d\lambda ^2 = 0 $
for all $w \in \Delta$. This means $d\lambda ^2 = 0$ on $S_p$. Under the
conditions of the theorem from (7) it follows that the restriction of
$\nabla R$ onto $S_p$ is zero. Then the equality (28) implies $DK = 0$ on
$S_p$, i.e. $S_p$ is a locally symmetric submanifold of $M$.

\hfill {$\square$}

The last question to consider is a theorem of Schur's type for the pointwise
constant relative sectional curvature (9).
\begin{thm}\label{th 3}
Let $(M,g)$ be a directed Riemannian manifold of pointwise constant relative
sectional curvature (9) and totally umbilical scalar distribution. Then the
curvature function $k$ is constant on the integral submanifolds of the scalar
distribution iff $\eta $ is closed ($\xi $ is geodesic).
\end{thm}

{\it Proof.} Writing the equality (24) in local coordinates
$$ \nabla _{i}\eta _{j} = \lambda (g_{ij} - \eta _{i}\eta _{j}) +
\eta _{i}\theta _{j}$$
we find
$$ \nabla _{i}\tau _{j} = \Vert d\tau \Vert _{i}\eta _{j} +
\Vert d\tau \Vert \lambda (g_{ij} - \eta _{i}\eta _{j}) -
\Vert d\tau \Vert \eta _{i}\theta _{j},$$
$$(dk + k \theta )\wedge \eta = 0.$$

The last equality shows that $d\ln k + \theta = 0$ on the integral
submanifolds $S_{p}$ of $\Delta $. Hence, $k$ is constant on $S_{p}$ iff
$\theta = 0$, i.e. $d\eta = 0$.

Finally the equalities
$$ g(\nabla _{\xi }\xi ,x) + \eta (\nabla _{\xi }x) = 0;$$
$$[\xi ,x] = \nabla _{\xi }x - \lambda x;$$
$$ \theta (x) = - \eta (\nabla _{\xi }x)$$
for all $x \in \Delta$ imply the condition $\theta = 0$ is equivalent to
the condition $\nabla _{\xi }\xi = 0,$ i.e. $\xi $ being geodesic.
\hfill {$\square$}

\begin{rem} In the given examples in section 3 a simple calculation
shows that
$$ \Vert d\tau \Vert ^2 = -{\frac{4(\tau - nB)^2 (\tau + 2B)}{(n-1)(n+2)}}
+ C(\tau - nB), \hspace{0,5cm} B, C = const.$$
Hence, the curvature function $k$ is a constant on $S_{p}$, but it is not
a global constant on $M$. Therefore Theorem 4.2 cannot be improved in this
direction.

It is interesting to find examples of directed Riemannian manifolds of
constant relative sectional curvature.
\end{rem}

\begin{rem} If $(M,g)$ is a surface with Gaussian curvature $K$ in
the Euclidean space, then its sectional 1-form $\varphi $ satisfies the
equality $\varphi = dK $ and consequently every surface is a directed
Riemannian manifold of pointwise constant relative sectional curvature
$k = \Vert dK\Vert $.

Hence, the surfaces of constant relative
sectional curvature are exactly the surfaces satisfying the condition
$\Vert grad \,K\Vert = const$.
\end{rem}
\vskip 4mm

The second author is partially supported by Sofia University Grant 99/2013.

\end{document}